\documentclass[12pt]{article}
\usepackage{latexsym, amsfonts, delarray}

\newtheorem{thm}{Theorem}[section]
\newtheorem{prop}[thm]{Proposition}
\newtheorem{cor}[thm]{Corollary}
\newtheorem{lem}[thm]{Lemma}
\newtheorem{defin}[thm]{Definition}
\newenvironment{unnumbered}[1]{\trivlist \item [\hskip \labelsep {\bf
#1}]\ignorespaces\it}{\endtrivlist}

\def \F {{\rm{I\hspace{-1mm}F}}}

\def \l {\left}
\def \r {\right}

\def \l {\left}
\def \r {\right}

\def\bigoh{\mathrm{O}}

\def\proof{\noindent \textsc{Proof:}\ }
\def\endproof{$\Box$ \medskip}

\def\Fq{\mathbb{F}_q}   

\def\Fqm{\mathbb{F}_{q^m}}   
\def\Fqr{\mathbb{F}_{q^r}}   
\def\Fqs{\mathbb{F}_{q^s}}
\def\Fqthree{\mathbb{F}_3}   
\def\Fqfive{\mathbb{F}_5}   
\def\Fqseven{\mathbb{F}_7}   
\def\FqXI{\mathbb{F}_{11}}   
\def\FqXIII{\mathbb{F}_{13}}

\def\ZZ{\mathbb{Z}}
\def\gal{\mathrm{Gal}}

\def\rand{\mathrm{and}}

\def\SL2{\mathrm{SL}_{2}}
\def\ZZ{\mathbb{Z}}

\def\proof{\noindent \textsc{Proof:}\ }
\def\endproof{$\Box$ \medskip}
\def\Pic{\mathrm{Pic}}

\begin{document}


\title{Coverings of curves  with asymptotically many rational points
\footnote{Research supported in part by the NSF grants
DMS96-22938 and DMS99-70651.}}

\author{\normalsize
Wen-Ching W Li
 and Hiren Maharaj\\ 
\small Department. of Mathematics \\[-5pt] \small Pennsylvania State
University \\[-5pt] 
\small University Park, PA 16802-6401 \\[-5pt] \small
wli@math.psu.edu, maharaj@math.psu.edu }

\date{25 August 1999}

\maketitle

\begin{abstract}
The number $A(q)$ is the upper limit of the ratio of the maximum 
number of points of a curve defined over $\Fq$  to the genus.
By constructing class field towers with good parameters
we present improvements of lower bounds of $A(q)$ for $q$ an odd power of a prime. 
\end{abstract}

\normalsize


\section{Introduction}
Given a finite field $\F_q$ of $q$ elements, by $K/\F_q$ we mean a global function field $K$ with full constant field $\F_q$, that
is, with $\F_q$ algebraically closed in $K$. A {\it rational place} of $K$
is a place of $K$ of degree 1. Write $N(K)$ for the number of rational
places of $K$ and $g(K)$ for the genus of $K$. According to the Weil-Serre
bound (see~\cite{ser1},~\cite{ser2}) we have
\begin{equation} \label{weil}
N(K)\le q+1+g(K)\lfloor 2q^{1/2}\rfloor ,
\end{equation}
where $\lfloor t\rfloor$ is the greatest integer not exceeding the real number
$t$.
 \begin{defin}
 For any prime power $q$ and any integer $g\ge 0$ put
\[
N_{q}(g)=\max N(K),
\]
where the maximum is extended over all global function fields $K$ of genus $g$ with full
constant field $\F_q$.
\end{defin}
In other words, $N_{q}(g)$ is the maximum number of
$\F_{q}$-rational points that a smooth, projective, absolutely irreducible
algebraic curve over $\F_q$ of genus $g$ can have. The following quantity was
introduced by Ihara~\cite{iha}.
\begin{defin}
For any prime power $q$ let
\[
A(q)=\lim \sup \nolimits_{g\to\infty}\frac{N_{q}(g)}{g}.
\]
\end{defin}
It follows from~(\ref{weil}) that $A(q)\le\lfloor 2q^{1/2}\rfloor$.
 Furthermore,
Ihara~\cite{iha}
showed that $A(q)\ge q^{1/2}-1$ if $q$ is a square. In the special cases
$q=p^2$ and $q=p^4$, this lower bound was also proved by Tsfasman,
Vl{\v a}dut, and Zink~\cite{tvz}. Hereafter, $p$ always denotes a 
prime number. Vl{\v a}dut and Drinfel'd~\cite{vd} established the
bound 
\begin{equation}\label{ub}
A(q)\le q^{1/2}-1  
\end{equation}
for all $q$. In particular this yields that
$A(q) = q^{1/2}-1$  if $q$ is a square.
Garcia and Stichtenoth~\cite{gs1},~\cite{gs3}  
proved that if $q$ is a square, then
$A(q)=q^{1/2}-1$ can be achieved by an explicitly constructed tower of
global function fields.
 
In the case where $q$ is not a square, no exact values of $A(q)$ are known,
but lower bounds are available which complement the general upper bound ~(\ref{ub}).
According to a result of Serre~\cite{ser1},~\cite{ser2} (see also~\cite{nx2})
   based on class field towers, we have
\begin{equation}\label{sb}
A(q)\ge c \log q
\end{equation}
with an absolute constant $c>0$. Zink~\cite{z} gave the best known lower bound for $p^3$:
\begin{equation}
A(p^{3})\ge\frac{2(p^{2}-1)}{p+2}. \nonumber
\end{equation}
Later, Perret~\cite{per} proved that if $l$ is a prime and if $q>4l+1$ and $q\equiv
1$ mod $l$, then
\begin{equation} \label{perb}
A(q^{l})\ge\frac{l^{1/2}(q-1)^{1/2}-2l}{l-1}.
\end{equation}
Niederreiter and Xing \cite{nx1} generalised and improved~(\ref{perb})
by establishing the following bounds.
If q is odd and $m\ge 3$  is an integer, then
\begin{equation} \label{newlb1}
A(q^{m})\ge\frac{2q}{\lceil 2(2q+1)^{1/2}\rceil+1}.
\end{equation}
 If $q\ge 4$  is even and $m\ge 3$  is an odd
integer, then
\begin{equation} \label{newlb2}
A(q^{m})\ge\frac{q+1}{\lceil 2(2q+2)^{1/2}\rceil+2}.
\end{equation}
As a consequence, they improved the 
Gilbert-Varshamov bound for sufficiently large composite nonsquare $q$ on a 
certain interval. 
Furthermore in~\cite{nx1}
they showed that $A(2)\ge\frac{81}{317}= 0.2555 \ldots$, 
$A(3)\ge\frac{62}{163}= 0.3803 \ldots$ and 
$A(5)\ge\frac{2}{3}= 0.666  \ldots$.

Denote the number of places of degree $r$  in a function field 
$F$ by $B_r(F)$ or simply $B_r$ if there is no danger of confusion.
Niederreiter and Xing  further extended their bounds~(\ref{newlb1}) and~(\ref{newlb2}) to
the following result in \cite{nx2}.
\begin{thm}\label{gennx}
Let $F/ \Fq $ be a global function field with $N \ge 1$ rational places. 
Let $r \ge 3 $ be an integer.
Suppose that  the ratio of class numbers
$h(F\Fqr)/h(F)$ is odd.
\newline $(1)$  If $q$ is odd and $B_r (F)\ge 2 (2N-1)^{1/2} + 3$, then 
\begin{equation}\label{newlb3}
A(q^r) \ge \frac{2(N-1)}{2g(F) + \lceil 2(2N-1)^{1/2} \rceil +1 }.
\end{equation}
$(2)$  If $q$ is even  and  $B_r (F)\ge 2 (2N-2)^{1/2} + 3$, then 
\begin{equation}\label{newlb4}
A(q^r) \ge \frac{N-1}{g(F) + \lceil 2(2N-2)^{1/2} \rceil +2}.
\end{equation}
\end{thm}
The bounds~(\ref{newlb1}) and~(\ref{newlb2}) follow from
the above theorem by considering the rational function 
field over $\Fq$. Using this theorem, they also found improved 
lower bounds for $A(q^3)$:
\begin{cor}\label{cubes}
$(1)$ If $q$ is a power of an odd prime $p$ and $p$ does not divide $\lfloor 2q^{1/2} \rfloor$, then 
\begin{equation} \label{cube1}
A(q^3) \ge \frac{2q + 4 \lfloor q^{1/2} \rfloor}{3 + \lceil 2(2q+4\lfloor q^{1/2} \rfloor +1)^{1/2} \rceil}.
\end{equation}
If $q$ is odd and $p$ divides  $\lfloor 2q^{1/2} \rfloor$, then 
\begin{equation} \label{cube2}
A(q^3) \ge \frac{2q + 4 \lfloor q^{1/2} \rfloor-4}{3 + \lceil 2(2q+4\lfloor q^{1/2} \rfloor -3)^{1/2} \rceil}.
\end{equation}
$(2)$ If $q\ge 4$ is even and $\lfloor 2q^{1/2} \rfloor$ is odd,  then 
\begin{equation} \label{cube3}
A(q^3) \ge \frac{q +  \lfloor 2q^{1/2} \rfloor}{3 + \lceil 2(2q+2\lfloor 2q^{1/2} \rfloor )^{1/2} \rceil}.
\end{equation}
 If $q\ge 4$ is even and $\lfloor 2q^{1/2} \rfloor$ is even,  then 
\begin{equation} \label{cube4}
A(q^3) \ge \frac{q +  \lfloor 2q^{1/2} \rfloor-1}{3 + \lceil 2(2q+2\lfloor 2q^{1/2} \rfloor -2)^{1/2} \rceil}.
\end{equation}
\end{cor}
A number $q$ is called special if $p$ divides  $\lfloor 2q^{1/2}\rfloor$ or  $q$ 
can be represented in one of the forms $n^2+1$, $n^2+n+1$, $n^2+n+2$ for some integer $n$.
They also proved that if $q\ge 11$ is odd, $\lfloor 2q^{1/2}\rfloor$  is even, and $q$
is not special, then  
\begin{equation} \label{cube5}
A(q^3) \ge \frac{2q + 4 \lfloor q^{1/2} \rfloor}{5 + \lceil 2(2q+4\lfloor q^{1/2} \rfloor +1)^{1/2} \rceil}.
\end{equation}

Recently Temkine~\cite{kem} extended Serre's lower bound ~(\ref{sb}) to

\begin{thm}\label{kemthm}
There exists an effective constant $c$ such that 
\begin{equation}\label{kembd}
A(q^r) \ge c r^2 \log q\frac{\log q}{\log r + \log q}.
\end{equation}
\end{thm}
It is also shown in~\cite{kem} that 
$A(3)\ge\frac{8}{17}= 0.4705 \ldots$ and 
$A(5)\ge\frac{8}{11}= 0.7272  \ldots$, thus improving the
corresponding bounds given in~\cite{nx1}. 

In this paper we employ class field towers to 
improve aforementioned lower bounds for $A(q)$ and to compute $A(p)$ for small primes $p$.  We also give an alternative  proof of Theorem~\ref{kemthm} with an explicit and improved constant $c$.
Finally, we present a lower bound for the $l$-rank of the $S$-divisor class group, similar to
the corresponding result from~\cite{nx1}.
More precisely, our results are as follows.

By using the explicit construction of ray class  fields of function fields via rank one 
Drinfeld modules  we prove the following generalisation of Theorem~\ref{gennx}.
\begin{unnumbered}{Theorems~\ref{qodd}}
Let $q$ be an odd prime power. Let $r$ be an odd integer at least $3$ and $s$ be a positive integer relatively
prime to $r$. Let $F/ \Fq $ be a global 
function field and let $N$ be the largest integer with the property 
that $B_s \ge N$ and $B_r $ $>$
$\left\lfloor (3+\lceil 2 (2N+1)^{1/2}\rceil )/ (r-2) \right\rfloor$.
Further suppose that $h(F\Fqr)/h(F)$ is odd.
Then we have
\begin{equation}
A(q^{rs}) \ge  \frac{4Ns}{4g(F) + 
      \left\lfloor \frac{3+\lceil 2 (2N+1)^{1/2}\rceil}{r-2} \right\rfloor 
                   + \lceil 2 (2N+1)^{1/2} \rceil} .
\end{equation}
\end{unnumbered}
By $\mathrm{f}(q)= \bigoh (\mathrm{g}(q))$ we mean that 
there is a constant $M>0$ such that  
$|\mathrm{f}(q)| \le M |\mathrm{g}(q)|$ for all sufficiently large $q$.
Immediate consequences of the above theorem are the following two corollaries.
\begin{unnumbered}{Corollary~\ref{qoddcor1}}
Let $q$ be an odd prime power. Let $r$ be an odd integer at least $3$ and $s$ be a positive integer relatively prime to $r$.  Let $F$ be the rational function field $\Fq (x) $ and suppose that 
\[
B_r >
\left\lfloor (3+\lceil 2 (2B_s+1)^{1/2}\rceil )/ (r-2) \right\rfloor.
\]
Then we have
\begin{equation}
A(q^{rs}) \ge  \frac{\sqrt{2}(r-2)}{r-1} \sqrt{s} q^{s/2} + \bigoh (1).
\end{equation}
\end{unnumbered}
For $r < s < 2r$  the conditions of this corollary are satisfied for 
all $q$ sufficiently large and the bound above improves
the bound~(\ref{newlb1}), which gives   
$A(q^{rs}) \ge  \frac{\sqrt{2}}{2} q^{s/2} + \bigoh (1)$.

Taking $F$ to be the rational function field and $s=1$,  one gets the following 
improvement of ~(\ref{newlb1}) for $r \ge 5$.
\begin{unnumbered}{Corollary~\ref{qoddcor} }
Let $q$ be an odd prime power.  Then for any odd integer $r \ge 3$ we have
\begin{equation}
A(q^r) \ge \frac{4q+4}
{\left\lfloor \frac{3+\lceil 2 (2q+2)^{1/2}\rceil}{r-2} \right\rfloor 
                   + \lceil 2 (2q+3)^{1/2} \rceil}.
\end{equation}
\end{unnumbered}
A better lower bound for $A(q)$ with $q$ even is derived from the 
following  generalisation of  Theorem~\ref{gennx}.
\begin{unnumbered}{Theorem~\ref{anyp}}
Let $F/ \Fq $ be a global function field of characteristic $p$. Let $r$ be an odd integer at least $3$ and $s$ be a positive integer relatively prime to $r$. Let $N$ be the largest integer such
that $B_s  \ge N$ and 
$B_r$ $ >$ $ \left\lfloor  \frac{6 + 2\lceil 2 \sqrt{2pN} \rceil}{r-1} \right\rfloor$.
If  $h(F\Fqr)/h(F)$ is not divisible by $p$, then 
\begin{equation}
A(q^r) \ge \frac{pNs}{p g(F) -p + 2(p-1)\left (3 + \l \lceil 2 \sqrt{pN} \r \rceil \right )}.
\end{equation}
\end{unnumbered}
\begin{unnumbered}{Corollary~\ref{qevencor}}
Let $q$ be a power of $2$. Let $r$ be an odd integer at least $3$ and $s$ be a positive integer relatively prime to $r$. Let $F$ be the rational function field $\Fq (x)$. Suppose that 
\[
B_r(F) > \left\lfloor \frac{6 + 2 \l\lceil 4 \sqrt{B_s(F)} \r\rceil}{r-1} \right\rfloor .
\]
Then we have
\begin{equation}
A(q^{rs}) \ge  \frac{\sqrt{2}}{4} \sqrt{s} q^{s/2} + \bigoh (1).
\end{equation}
\end{unnumbered}
For $r < s < 2r$  the conditions of this corollary are satisfied for 
all even $q$ sufficiently large and the bound above improves
the bound~(\ref{newlb2}), which gives   
$A(q^{rs}) \ge  \frac{\sqrt{2}}{4} q^{s/2} + \bigoh (1)$.

By applying Theorem~\ref{anyp} to Deligne-Lusztig curves
in characteristic 2,
we obtain the following bound which improves ~(\ref{newlb2}), 
(\ref{cube3}) and (\ref{cube4}).
\begin{unnumbered}{Theorem~\ref{qeven}}
Let $q$ be a power of $2$.
For $r \ge 5$ odd and all $q$ sufficiently large we have
\begin{equation}
A(q^r)  \ge    \frac{2q^2}{\sqrt{2q}(q-1) + 2 \lceil 2\sqrt{2}q \rceil +4}.
\end{equation}
For $r = 3$ and all $q$ sufficiently large we have
\begin{equation}
A(q^3)   \ge  \frac{2q^2}{\sqrt{2q}(q-4) + 8 \lceil \sqrt{2}q \rceil +16}.
\end{equation}
\end{unnumbered}
The same ideas involved in the proof of the lower bound of $A(q^3)$ for 
$q$ even can be used to prove the following bounds  which improve the bounds of 
Corollary~\ref{cubes} and the bound~(\ref{cube5}) for 
characteristics 3, 5, and 7.
\begin{unnumbered}{Theorem~\ref{357}}
Let $q$ be a power of $p=$ $3$, $5$ or $7$. Then for all $q$ sufficiently large we have
\begin{equation}   
A(q^3)   \ge  \frac{2(q^2+p^2)}{\sqrt{pq}(q-p^2) + 4p(p-1) \lceil \sqrt{q^2/p + p}\rceil +10p^2-12p}.
\end{equation}
\end{unnumbered}
All of the above lower bounds for $A(q^r)$ are good for large $q$. The next result, which
is a generalisation of the bound~(\ref{sb}),  is better for $r$ large.
Theorem~\ref{kemthm} is a consequence of this.
\begin{unnumbered}{Theorems~\ref{kemqodd} and~\ref{kemqeven}}
Let $0 < \theta < 1/2$.
Then for all sufficiently large odd $q^r$  we have 
\begin{equation}\label{kem1}
A(q^r) \ge
\frac{((\lfloor \theta r   \log q \rfloor-3)^2 - 4)r}
{2(\lceil 2\log r/ \log q \rceil +1)(\lfloor \theta r   \log q \rfloor +1)-6}.
\end{equation}

\noindent
For all sufficiently large  even $q^r$ we have 
\begin{equation}\label{kem2}
A(q^r) \ge
\frac{\l ( \lfloor \theta r \log q \rfloor-2 \r )^2r}
{4(\lceil 2\log r/ \log q \rceil +1)(\lfloor \theta r \log q \rfloor +1)-8}.
\end{equation}
\end{unnumbered}

Using Tate cohomology, Niederreiter and Xing obtained  
a lower bound for the $l$-rank of the 
$S$-divisor class group $Cl_S$ (see Proposition~\ref{nxlb}).
In section~\ref{smallp} we present a proof of the following similar result. 
\begin{unnumbered}{Proposition~\ref{nxlb2}}
Let $F/ \Fq$ be a global function field and $K$  a finite abelian
extension of $F$. Let $T$ be a finite nonempty set of places of $F$ and $S$ the set of
places of $K$ lying over those in $T$. 
If at least one place in 
$T$ splits completely in $K$, then  for any prime $l$ we have
\begin{equation}
d_{l}Cl_{S}\ge 
       \sum_{P}d_{l}G_{P}-(|T|-1+d_{l}\mathbb{F}_{q}^*   )-d_{l}G,
\end{equation}
 where  $G= \gal (K/F)$, $G_P$ is the inertia subgroup at the
place $P$ of $F$, and $d_{l}X$ denotes the $l$-rank of an abelian group $X$. 
The sum is extended over all places $P$ of $F$.
\end{unnumbered}
It is easily shown that this lower bound coincides with that of Niederreiter and Xing. 
 The proof of the bound, which uses
narrow ray class fields,  reveals that the lower bound 
is really a lower bound of the $l$-rank of $\gal(K'/K)$, where $K'$ is 
the maximal subfield of the $S$-Hilbert class field of $K$ which is an abelian extension of $F$. 
Finally, in section~\ref{smallp}, lower bounds for $A(p)$ for small primes $p$ 
are computed.
We obtain
$A(7) \ge 9/10$,
$A(11) \ge 12/11 = 1.0909 \ldots$ and 
$A(13) \ge 4/3$ and $A(17) \ge 8/5$.


\section{Background on class field theory}
\subsection{Hilbert class fields}\label{hcf}
We recall, without proof, 
 some basic facts about Hilbert class fields. The reader is referred to \cite{rosen}
for more details.
Let $K/\F_q$ be a global function field with full constant field $\F_q$. 
Let $S$ be a finite nonempty set of places of $K$ and $O_S$  the $S$-{\it 
integral ring} of $K$, i.e., $O_S$ consists of
all elements of $K$ that have no poles outside $S$. Denote by $O_{S}^*$ 
the group of units in $O_S$. If $S$ consists of just one element $P$, then
we write $O_P$ and $O_P^*$ for  $O_S$ and $O_S^*$.
The $S$-{\it Hilbert class field} $H_S$ of $K$
 is the maximal unramified  abelian extension of $K$
(in a fixed separable closure of $K$) in which all places in $S$ split
completely. 
The galois group of $H_S/K$, denoted by $Cl_S$, is isomorphic to the class
group of $O_S$ (see~\cite{rosen}); its order is the class number  $h(O_S)$.
If $S=\{P\}$ consists of one element, then $h(O_S) = d h(K)$ with 
$d = \deg P$ and $h(K)$ the divisor class number of $K$.

Now we define the $S$-{\it class field tower} of $K$. Let $K_1$ be the
$S$-Hilbert class field $H_S$ of $K$ and $S_1$ the set of places of $K_1$
lying over those in $S$. Recursively, we define $K_i$ to be the 
$S_{i-1}$-Hilbert class field of $K_{i-1}$ for $i\ge 2$ and $S_i$ to be the
set of places of $K_i$ lying over those in $S_{i-1}$. Then we get the 
$S$-{\it class field tower} of $K$:
$K=K_{0}\subseteq K_{1}\subseteq K_{2}\subseteq \ldots.$
The tower is infinite if $K_{i}\not=K_{i-1}$ for
all $i\ge 1$. The following proposition, known to Serre and proved by Schoof in ~\cite{schoof}, provides a sufficient condition for a class field tower to be infinite. It is a basic tool for our 
work, and it leads to the stated lower bound for $A(q)$. 
For a prime $l$ and an
abelian group $B$, denote by $d_{l}B$ the $l$-rank of $B$.
\begin{prop}~\cite{schoof} \label{golshaf}  
Let $K/\F_q$ be a global function field of genus $g(K)>1$ and let $S$ be a nonempty set of
 places of $K$. Suppose that there exists a prime $l$ such that
\begin{equation} \label{gscond}
 d_{l}Cl_{S}\ge 2+2\l(d_{l}O_{S}^{*}+1\r)^{1/2}.
\end{equation}
Then $K$ has an infinite $S$-class field tower. Furthermore if $S$ consists of 
only rational places, then
\[ 
A(q)\ge\frac{|S|}{g(K)-1}.
\]
\end{prop}
The $l$-rank of $O_{S}^*$ can be determined. Dirichlet's unit theorem asserts that 
$O_{S}^{*}\simeq\F_{q}^{*}\times{\bf Z}^{|S|-1}$,
and therefore
$$d_{l}O_{S}^{*}= \left\{ \begin{array}{l@{\quad}l}
|S| & \hbox{if} \; l|(q-1),\\
|S|-1 & \hbox{otherwise.} \end{array} \right.$$

\subsection{Narrow ray class fields} \label{nrcf}

Since the explicit constructions of ray class fields
 via Drinfeld modules of rank 1 will be used, we recall
the results and basic definitions. 
For more information the reader may consult ~\cite{hay},~\cite{rosen} 
and~\cite{goss}.
We shall follow the same notation as~\cite{nx3}.
 
Let $F/\Fq$ be a global function field. We distinguish a place
$\infty$ of $F$ and let $A$ be the subring of $F$ consisting of all
the functions which are regular away from $\infty$. Then the 
Hilbert class field $H_A$ of $F$ with respect to $A$ is the
maximal unramified abelian extension of $F$ (in a fixed separable
closure of $F$) in which the place 
$\infty$ splits completely. The galois group of $H_A/F$ is 
isomorphic to the fractional ideal class group $\Pic(A)$ 
of $A$. If the degree of $\infty$ is 1 then the degree $[H_A:F]$
is the divisor class number $h(F)$ of $F$.

We fix a sign function sgn and let $\phi$ be a rank 1 sgn normalised 
Drinfeld $A$-module over $H_A$. The additive group of the algebraic 
closure $\overline{H_A}$ of $H_A$ forms an $A$-module under the 
action of $\phi$. For any nonzero integral ideal $M$ of $A$, the $M$-torsion module $\Lambda (M) = \{ u \in \overline{H_A} : 
\phi_M(u)= 0 \}$ is a cyclic $A$-module
which is isomorphic to the $A$-module $A/M$ and has $\Phi_q (M) := |(A/M)^*|$
generators, where $(A/M)^*$ is the group of units of the ring $A/M$.
Let $\mathcal{I}(A)$ be the fractional ideal group of $A$ and let 
$\mathcal{I}_M(A)$ be the subgroup of fractional ideals in $\mathcal{I}(A)$ 
prime to $M$. Define the quotient group $\Pic_M(A) = \mathcal{I}_M(A)/
\mathcal{R}_M(A)$, where $\mathcal{R}_M(A)$ is the subgroup consisting of 
all principal ideals $bA$ with $\mathrm{sgn}(b)=1$ and $b \equiv 1 \hbox{ mod } M$.

The field $F_M = H_A( \Lambda (M))$ generated by the elements of $\Lambda (M)$
over $H_A$ is called the narrow ray class field of $F$ modulo $M$.
The extension $F_M$ is unramified away from $\infty$ and  the 
prime ideals in $A$ which divide $M$.
In fact, the  maximal 
subfield  in which $\infty$ splits completely is the ray class
field of $F$ with conductor $M$.
We summarize below the main results from~\cite{hay} which will be used later in the paper.
\begin{prop}
Let $F_M = H_A( \Lambda (M))$ be the narrow ray class field of $F$ modulo
$M$. Then:

\noindent
$1.$ $F_M/F$ is an abelian extension and there is an isomorphism 
$\sigma: \mathrm{Pic}_M(A) \to \gal (F_M/F)$ given by
$\sigma_I \phi = I * \phi$ for any ideal $I$ in $A$ prime to $M$, 
and $ \lambda ^{\sigma_I} = \phi_I(\lambda)$ for any generator
$\lambda$ of the cyclic $A$-module $\Lambda (M)$. Moreover for any 
ideal $I$ in $A$ prime to $M$, the corresponding Artin automorphism 
of $F_M/F$ is $\sigma_I$.  \newline
$2.$ the muliplicative group $(A/M)^*$ is isomorphic to $\gal (F_M/ H_A )$
by means of $b \mapsto \sigma_{bA}$, where $b \in A$ is prime to $M$
and satisfies $\mathrm{sgn}(b) = 1$.  \newline
$3.$ both the decomposition group and the inertia group of 
$F_M/F$ at $\infty$ are isomorphic to 
the multiplicative group $\mathbb{F}_{q}^*$. \newline
$4.$ if $M$ decomposes into the product $M= \prod_{i=1}^t P_i^{e_i}$ 
of distinct prime ideals in $A$ with $e_i \ge 1$,  then 
$F_M$ is the composite of the fields  $H_A( \Lambda (P_1^{e_1}))$,  
$H_A( \Lambda (P_1^{e_2}))$, ..., $H_A( \Lambda (P_t^{e_t}))$. The 
order of the inertia group of $F_M/F$ at $P_i$ is 
$\Phi_q (P_i^{e_i})$, $i = 1,2, \ldots, t$.
\end{prop} 


\section{General lower bounds for $A(q^r)$}

Recall that $B_r(F)$ or simply $B_r$ denotes the number of places of degree $r$  in a function field 
$F$.  
The following estimate of the size of $B_r$ was proved in \cite{stich} (Corollary~V.2.10).
\begin{prop}\label{brrrprop}
For a global function field 
$F/ \Fq$ we have
\begin{equation}\label{brrr}
\l | B_r - \frac{q^r}{r} \r |  \le \left (\frac{q}{q-1} + 2g(F)\frac{q^{1/2}}{q^{1/2}-1} \right ) \frac{q^{r/2}-1}{r}
    <   (2+7g(F)) \frac{q^{r/2}}{r},
\end{equation}
where $g(F)$ is the genus of $F$.
\end{prop}
By $\mathrm{f}(q)= \bigoh (\mathrm{g}(q))$ we mean that 
there is a constant $M>0$ such that  
$|\mathrm{f}(q)| \le M |\mathrm{g}(q)|$ for all sufficiently large $q$.
Thus Proposition~\ref{brrrprop} implies that $B_r(\Fq (x))  =  q^r/r + \bigoh (q^{r/2})$.

\subsection
{Lower bounds for $A(q^r)$ with $q$ odd}\label{nlbo}

We start by proving a general theorem from which several improvements on lower bounds will be derived. Narrow ray class fields are used to construct infinite
towers of function fields over $\Fqr$. 

For an odd integer $r \ge 3 $ denote by $F_r$ the extension of $F$ by the constant field $\Fqr$. Let 
 $s$ be a positive integer relatively prime to $r$.
Then all the places of degree $s$ in  $F$ can be viewed naturally as degree $s$ places of $F_r$. 
As in section 2.2, let $A_r$ be the  subring of $F_r$ consisting of 
elements which are regular outside a chosen place $\infty$.
Any place $Q$ of degree $r$ decomposes
into a product of $r$ prime ideals of degree 1 in $A_r$.
In case no confusion can arise, we denote both the ideal
$Q \cdot A$ and $Q \cdot A_r$ simply by $Q$. Thus
$(A_r/Q)^* = (A_r / Q\cdot A_r)^*$ can be regarded as a 
subgroup of $\Pic_Q (A_r) = \Pic_{Q \cdot A_r}(A_r)$ in a canonical
way. 

As described in~\cite{nx3}, the group $\Pic_Q(A)$ can also
be viewed as a subgroup of $\Pic_Q(A_r)$ in a natural way.
This is  explained in the language of algebraic curves 
as follows. Let $C$ be an algebraic curve over $\Fq$ with 
function field $F$. If we view $C$ as a curve over $\mathbb{\overline{F}}_q$,
then a divisor $D$ on $C/\mathbb{\overline{F}}_q$ is a divisor of $F$
if and only if $D$ is $\Fq$-rational, that is, invariant under the action of $\gal (\mathbb{\overline{F}}_q / \Fq)$. 
Hence $\Pic_Q(A)$ can be described as the group of all 
$\Fq$-rational divisors on $C/\mathbb{\overline{F}}_q$
prime to $Q$ and $\infty$, from which we factor out the group of all divisors
$(c)_0$ with $c \in F$, $\mathrm{sgn}(c)=1$, and $c \equiv 1 \mathrm{\; mod \; } Q$ where
$(c)_0$ is the divisor corresponding to the principal ideal $cA$.
Since $\Pic_Q(A_r)$ has a similar description it follows that 
$\Pic_Q(A)$  is a subgroup of $\Pic_Q(A_r)$ in a natural way.

We will use the following result from~\cite{nx2}.
\begin{lem}\label{forqodd}
Given a place $Q$ of degree $r$ of $F$, let $E_r = H_{A_r}(\Lambda(Q\cdot A_r))$ be the narrow
ray class field of $F_r$ modulo $Q \cdot A_r$. Let $L$ be the subfield
of $E_r/F_r$ fixed by the subgroup $\Pic_Q(A)$ of 
$\gal (E_r/F_r)$, and let $K/F_r$ be the maximal unramified
extension of $F_r$ in $L$. Then the degree of the 
extension $K/F_r$ is $h(F_r)/h(F)$.
\end{lem}

Let $L$ be as in this lemma. 
Observe that a degree $s$ place of $F_r/\Fqr$ different from $\infty$
splits completely in $L/F_r$ if and only if its Artin automorphism
is contained in $\Pic_Q(A)$, and this happens if and only if 
the restriction of this place to $F/ \Fq$ is a place of degree $s$. This fact will be used repeatedly. 
\begin{thm}\label{qodd}
Let $q$ be an odd prime power. Let $r$ be an odd integer at least $3$ and $s$ be a positive integer relatively prime to $r$. Let $F/ \Fq $ be a global 
function field and let $N$ be the largest integer such 
that $B_s \ge N$ and $B_r $ $>$
$\left\lfloor (3+\lceil 2 (2N+1)^{1/2}\rceil )/ (r-2) \right\rfloor$.
Further suppose that $h(F_r)/h(F)$ is odd.
Then we have
\begin{equation}
A(q^{rs}) \ge  \frac{4Ns}{4g(F) + 
      \left\lfloor \frac{3+\lceil 2 (2N+1)^{1/2}\rceil}{r-2} \right\rfloor 
                   + \lceil 2 (2N+1)^{1/2} \rceil} .
\end{equation}
\end{thm}
\proof 
Put $ n = 
\left\lfloor (3+\lceil 2 (2N+1)^{1/2}\rceil )/ (r-2) \right\rfloor$ and 
let $Q_{1},\ldots,Q_{n+1}$ be $n+1$ distinct places of degree
$r$ in $F$.
 Then each $Q_i$ decomposes into a product
$Q_{i}A_r=\prod_{j=1}^{r}Q_{ij}$
of $r$ distinct prime ideals of degree one in $A_r$. 
For each $i$ consider the narrow ray class field
$E_{r}^{(i)} = H_{A_r}(\Lambda (Q_i))$ of $F_r$ modulo $Q_i$.
We use the abbreviations $h=h(F)$ and $h_r = h(F_r)$ for the remainder of the proof.

Let $I_i$ be the inertia group of $E_{r}^{(i)}/F_r$ at
$\infty$  and $L_i$ the
subfield of $E_{r}^{(i)}/F_r$ fixed by the subgroup 
$I_{i}\cdot \Pic_{Q_i} (A)$
of Gal$(E_{r}^{(i)}/F_r)$.
Since $|\Pic_{Q_i} (A)|= h(q^r-1)$, $|I_i|= q^r-1$ and 
$| I_{i}\cap \Pic_{Q_i} (A)| =q-1$, it follows that
$| I_{i}\cdot \Pic_{Q_i} (A)| =h(q^r-1)^2/(q-1)$.
Hence $[L_i:F_r] = (h_r/h)(q-1)(q^r-1)^{r-2}$.
The order of the inertia group
of $Q_{ij}$ in $L_{i}/F_r$ divides $|(A_r/Q_{ij})^*|=q^{r}-1$ 
for each $1\le j\le r$, and therefore
the inertia groups of $Q_{i3},\ldots,Q_{ir}$  in $L_{i}/F_r$  generate a subgroup $G_i$ of
$\gal (L_{i}/F_r)$ of order dividing $(q^{r}-1)^{r-2}$. Let $J_i$ be the subfield
of $L_{i}/F_r$ fixed by $G_i$, then $(h_r/h)(q-1)$ divides 
the degree of the extension $J_{i}/F_r$. The only
possible ramified places in $L_{i}/F_r$ are $Q_{i1}$ and  $Q_{i2}$.
 
Let $K_{i \, 2}$ be a quadratic extension of $F_r$ in $L_i$. (The reason for our choice of notation will become clear.) 
The only possible
ramified places in $K_{i \, 2}/F_r$ are $Q_{i1}$ and  $Q_{i2}$.
On the other hand, since $h_r/h$ is odd, by Lemma~\ref{forqodd}, the field
$K_{i \, 2}$ is not contained in the maximal unramified extension of
$F_r$ in the subfield of $E_r^{(i)}$ fixed by $\Pic_{Q_i} (A)$. In other words,
$  K_{i \, 2}/F_r$ is ramified. Thus at least one of the places
$Q_{i1}$, $Q_{i2}$ is ramified in $  K_{i \, 2}/F_r$.
It is impossible for exactly one of these places to ramify in 
$  K_{i \, 2}/F_r$, for otherwise the Hurwitz genus formula would
yield
$2g(K_{i \, 2})-2 = 2(2g(F_r)-2) + (2-1)\cdot 1,$
contradicting the integrality of $g(K_{i \, 2})$.
Since $[K_{i \, 2}:F_r]=2$ and $q^r$ is odd, $K_{i \, 2}/F_r$
is a Kummer extension and we can write $K_{i \, 2} = F_r(y_{i \, 2}),$
where $y_{i \, 2}^2$ equals an element $ u_{i \, 2} \in F_r$.
Since $Q_{i1}$ and  $Q_{i2}$ are the only places of $F_r$ that ramify
in $K_{i \, 2}$, it follows that for each place $P$ of $F_r$, the valuation at $P$ of $u_{i \, 2}$, noted $v_P(u_{i \, j})$, is odd
only when $P= Q_{i1}, \, Q_{i2}$.

Thus, by the above argument, for $1 \le i \le B_r$ and $2 \le j \le r-1 $
we can form extensions 
$K_{i \, j} = F_r(y_ {i \, j}),$ where 
\[
y_ {i \, j}^2 = u_{i \, j}
\]
and in $K_{i \, j}/F_r$ the only two places that ramify are 
$Q_{i \, 1}$ and $Q_{i \, j}$ so that the $u_{i \, j}$ have the 
property that
$v_P(u_{i \, j})$ is odd only when $P= Q_{i1}, \, Q_{ij}$.

Let $K_j'$ denote the compositum of the fields $K_{i \, j'}$ for $2 \le j' \le r-1$ and $j' \ne j$. Observe that $K_{i \,  j} \cap K_j' = F_r$ because the place $Q_{ij}$ is totally ramified in 
$K_{i \, j}/F_r$ but unramified in $K_j'/F_r$.

Before going further, we introduce more notation.
Put $A = 3 + \lceil 2 \sqrt{2N+1} \rceil$. 
Let $t_1$ be the integer $A - n(r-2)$. Thus $0 \le t_1 <  r-2$.
If $t_1=0$, set $t=0$; if $t_1>0$ and  $t_1$ is even, set 
$t=t_1$; otherwise set $t=t_1 +1$.
Define the sets
\[
Z = \{ (i,j) |\, 1 \le i \le n,  2 \le j \le r-1 \}  \cup \{(n+1, j) |\,  
2 \le j \le t \},
\]
where the second set is empty if $t = 0$ and
\[
Z' = Z \cup \{(i, 1) |\, (i, 2) \in Z \}.
\] 
Form the extensions
\[
K = F_r (\{  y_{i \, j}|  (i,j) \in Z \}). 
\]
and $L=F_r(y)$ with
\begin{equation}\label{cyclic}
y^2 =  \prod_{(i,j) \in Z } u_{i \, j}(x). 
\end{equation}
The galois group of $K/F_r$ is elementary abelian of exponent 2.

Since $y^2$ equals  a product of some of the $y_{i \, j}^2$'s, it follows that
$L$ is a subfield of $K$.  
Observe from construction that if $Z$ contains one pair $(i, j)$, then it contains an odd number of pairs with the first component $i$. Consequently the places $Q_{ij}$ with $(i, j) \in Z'$  are ramified 
in $K/F_r$ with ramification index 2 by repeated 
 application of Abhyankar's Lemma
(see~\cite{stich} chapter III).  The same happens in the 
extension $L/F_r$. Therefore the  extension $K/L$ is unramified. 

Let $T'$ be a set of $N$ places in $F$ of degree $s$ and $T$ the set
of $N$ places in $F_r$ which lie above those in $T'$. 
Then from the remarks preceding this theorem we know that all the places
in $T$ split completely in $L/F_r$.
Let $S'$ denote the $2N$ places in $L$ which lie above the 
places in $T$.
Now $K/L$ is an unramified abelian extension in which
the places in $S'$  split completely.
Moreover we have 
$d_{l}Cl_{S'} \ge d_{l}\gal (K/L) = (n(r-2)+ t)-1$.
Since
\[
n(r-2)+ t -1 \ge n(r-2)+ t_1 - 1 =  2 + 2 \sqrt{2N+1}=2 + 2 \sqrt{|S'|+1} ,
\]
it follows from  Proposition~\ref{golshaf} that
$L$ has an infinite $S'$-Hilbert class field tower.

By the Hurwitz genus formula we have,
\begin{eqnarray}
2g(L)-2  & =   & 2(2g(F_r)-2) + n(r-1)+t  \nonumber \\
         & =   & 4g(F) + n-4 + n(r-2)+t          \nonumber \\
         & \le & 4g(F) + n-4 + A+1         \nonumber \\
         & =   & 4g(F) + 
      \left\lfloor \frac{3+\lceil 2 (2N+1)^{1/2}\rceil}{r-2} \right\rfloor 
                   + \lceil 2 (2N+1)^{1/2} \rceil . \nonumber
\end{eqnarray}

Passing to the constant field extension $L\Fqs$, each place of degree
$s$ in $S'$ splits into $s$ places of degree 1 in $L\Fqs$  
and it is  easily seen that  $L\Fqs$
 has an infinite $S$-Hilbert class field tower where $S$ is the set 
of those places
of $L\Fqs$ which lie above those in $S'$. We thus get, again by 
Proposition~\ref{golshaf},
that
\begin{eqnarray}
A(q^{rs})& \ge & \frac{|S|}{g(L\Fqs)-1} = \frac{|S'|s}{g(L)-1}  \nonumber  \\ 
& \ge & \frac{4Ns}{4g(F) + 
      \left\lfloor \frac{3+\lceil 2 (2N+1)^{1/2}\rceil}{r-2} \right\rfloor 
                   + \lceil 2 (2N+1)^{1/2} \rceil}, \nonumber 
\end{eqnarray}
as desired.\endproof

When applied to rational function fields, the theorem above yields the following lower bounds. 

\begin{cor}\label{qoddcor1}
Let $q$ be an odd prime power. Let $r$ be an integer at least $3$ and $s$ be a positive integer relatively
prime to $r$. Let $F$ be the rational function field $\Fq (x) $. Suppose that 
\[
B_r >
\left\lfloor  \l (3+\lceil 2 (2B_s +1)^{1/2}\rceil \r )/ (r-2) \right\rfloor.
\]
Then we have
\begin{equation}\label{substa}
A(q^{rs}) \ge  \frac{\sqrt{2}(r-2)}{r-1} \sqrt{s} q^{s/2} + \bigoh (1).
\end{equation}
\end{cor}
For $r < s < 2r$  the conditions of Corollary~\ref{qoddcor1} are satisfied for 
all $q$ sufficiently large and the bound~(\ref{substa}) improves
the bound~(\ref{newlb1}) which gives 
$A(q^{rs}) \ge  \frac{\sqrt{2}}{2} q^{s/2} + \bigoh (1)$.

Taking $F$ to be the rational function field and $s=1$,  one gets the following 
bound which improves~(\ref{newlb1}) for $r \ge 5$.
\begin{cor}\label{qoddcor} 
Let $q$ be an odd prime power.  Then for any odd integer $r \ge 3$ we have
\begin{equation}\label{nx2}
A(q^r) \ge \frac{4q+4}
{\left\lfloor \frac{3+\lceil 2 (2q+2)^{1/2}\rceil}{r-2} \right\rfloor 
                   + \lceil 2 (2q+3)^{1/2} \rceil}.
\end{equation}
\end{cor}
We remark that in the case that $F$ is the rational function
field $\Fq (x)$, the function fields $K_{i j}$, and hence $L$, 
 of Theorem~\ref{qodd}
can be explicitly defined. See~\cite{wli} for the details.

%
%

By using similar ideas as in the proof of Theorem~\ref{qodd}
one can prove the next theorem. Instead of using the modulus
$Q$, we use $Q^2$. Moreover Artin-Schreier extensions are used instead of 
Kummer extensions.

\begin{thm}\label{anyp}
Let $F/ \Fq $ be a global 
function field of characteristic $p$. Let $r$ be an odd integer at
least $3$ and $s$ be a positive integer relatively prime to $r$. Let $N$
be the largest integer such that $B_s  \ge N$ and 
$B_r$ $ >$ $ \left\lfloor  \frac{6 + 2\lceil 2 \sqrt{2pN} \rceil}{r-1} \right\rfloor$.
If  $h(F\Fqr)/h(F)$ is not divisible by $p$, then 
\begin{equation}
A(q^r) \ge \frac{pNs}{p g(F) -p + 2(p-1)\left (3 + \l \lceil 2 \sqrt{pN} \r \rceil \right )}.
\end{equation}
\end{thm}
One obtains similar corollaries as before.
\subsection
{Lower bounds for $A(q^r)$ with $q$ even}\label{nlbe}
 
We start with a consequence of Theorem~\ref{anyp} for the case of even $q$.
\begin{cor}\label{qevencor}
Let $q$ be a power of $2$. Let $r$ be an odd integer at least $3$ and let $s$ be a positive integer relatively prime to $r$. Let $F$ be the rational function field $\Fq (x)$. Suppose that 
\[
B_r > \left\lfloor \frac{6 + 2 \l\lceil 4 \sqrt{B_s} \r\rceil}{r-1} \right\rfloor .
\]
Then we have
\begin{equation}\label{substa2}
A(q^{rs}) \ge  \frac{\sqrt{2}}{4} \sqrt{s} q^{s/2} + \bigoh (1).
\end{equation}
\end{cor}
For $r < s < 2r$  the conditions of Corollary~\ref{qevencor} are satisfied for 
all $q$ sufficiently large and the bound~(\ref{substa2}) improves
the bound~(\ref{newlb2}), which gives $A(q^{rs}) \ge  \frac{\sqrt{2}}{4} q^{s/2} + \bigoh (1)$.

Letting $s=1$,  we obtain the following 
result which is similar to the bound~(\ref{newlb4}) of Theorem~\ref{gennx}. 
\begin{thm}\label{qeven}
Let $q$ be a power of $2$. Let $F/ \Fq $ be a global 
function field with $N$ rational places.
Suppose that 
$B_r > \left\lfloor  \frac{6 + 2\lceil 2 \sqrt{2N} \rceil}{r-1} \right\rfloor$
and that  the ratio of class numbers
$h(F\Fqr)/h(F)$ is not divisible by $2$.
Then 
\begin{equation}
A(q^r) \ge \frac{N}{g(F) + \l \lceil 2 \sqrt{2N} \r \rceil  + 2}.
\end{equation}
\end{thm}
\indent Next we use this theorem to  prove
a lower bound for $A(q^r)$ which improves the bound~(\ref{newlb2}).

\begin{lem}\label{div}
Let $F/\Fq$ be a function field with at least one place of degree $r$ and more than one 
rational place. Then $h(F_r)/h(F)$ divides the class number
$h(O_S)$, where $S$ consists of all but one rational places in $F$ (viewed in $F_r$).
\end{lem}
\proof
Let $Q$ be a place of degree $r$ in $F$. Denote by $\infty$ the rational  place of $F$ not contained in $S$ and define the ring $A_r$ as before. Let $E_r = H_{A_r}(\Lambda(Q\cdot A_r))$ be the narrow
ray class field of $F_r$ modulo $Q \cdot A_r$. Let $L$ be the subfield
of $E_r/F_r$ fixed by the subgroup $\Pic_Q(A)$ of 
$\gal (E_r/F_r)$, and let $K/F_r$ be the maximal unramified
extension of $F_r$ in $L$. 
Then, from the remarks preceding Theorem~\ref{qodd},     
all places in $S$ split completely in $K/F_r$ and, 
from Lemma~\ref{forqodd},  the degree of the 
extension $K/F_r$ is $h(F_r)/h(F)$.
On the other hand the degree of the maximal unramified abelian extension 
of $F_r$  in which all the 
places in $S$ split completely is $h(O_S)$ (see section~\ref{hcf}). Hence $h(F_r)/h(F)$ divides 
$h(O_S)$.\endproof
 
\noindent
We will use the following  result proved by Rosen~\cite{rosen}.
\begin{prop}\label{rose}
Let $L/K$ be a galois extension with degree a power of a prime $l$. Let $S$ be a finite nonempty set of places of $K$. Suppose that every place in $S$ splits
completely in $L$ and that at most one place of $K$ ramifies
in $L$. If $S'$ is the set of primes of $L$ which lie above those in $S$, 
then $l|h(O_{S'})$ implies $l|h(O_{S})$.
\end{prop}

\begin{thm}\label{qevenagain}
Let $q$ be a power of $2$.
For $r \ge 5$ odd and $q$ sufficiently large we have
\[
A(q^r)  \ge   \frac{2q^2+2}{\sqrt{2q}(q-1) + 2 \lceil 2\sqrt{2q^2+2} \rceil +4}.
\]
For $r = 3$ and $q$ sufficiently large we have
\[    
A(q^3)   \ge  \frac{2q^2+8}{\sqrt{2q}(q-4) + 8 \lceil \sqrt{2q^2+8}\rceil +16}.
\]
\end{thm}
\proof
Set $K=\Fq(x)$ and $K_r = \Fqr (x)$.  Write $q=2^{2m+1}= 2q_0^2$ and define the extension 
$L=K(y)$ by 
\[
y^q+y = x^{q_0}(x^q+x).
\] 
Then $L$ has degree $q$ over $K$, it is totally ramified at $\infty$ and 
totally split at all other places of degree 1. Thus $L$ has $N=q^2+1$ rational 
places. As computed in ~\cite{gs4}, $L$ has genus $q_0(q-1)$. Let $L_r = K_r(y)$.

Let $r \ge 3$ be an odd number.
In order to apply Theorem~\ref{qeven} to the function field $L$, we must 
show that the number of places of degree $r$ in 
$L$ satisfies the condition 
$B_r(L)  >
\left\lfloor  (6 + 2\lceil 2 \sqrt{2N} \rceil)/(r-1) \right\rfloor
=\left\lfloor  (6 + 2\lceil 2 \sqrt{2q^2+2} \rceil)/(r-1) \right\rfloor$.
By Proposition~\ref{brrrprop}  
we have $B_r(L) > (1/r)(q^r -(7/\sqrt{2})q^{\frac{r+3}{2}}+(7/\sqrt{2})q^{\frac{r+1}{2}}
-2q^\frac{r}{2})$, hence the desired condition is satisfied for $r \ge 5$ and 
$q$ sufficiently large. 
The extension $L_r/K_r$ satisfies the conditions of Proposition~\ref{rose}
with  $l=2$. Since $h(K_r)=1$, it follows from Proposition~\ref{rose} and Lemma~\ref{div}
that $h(L_r)/h(L)$ is odd. 

Thus by Theorem~\ref{qeven},  for $r\ge 5$ and all $q$ sufficiently 
large  we get
\begin{eqnarray}
A(q^r)  & \ge  & \frac{N}{g(L) + \l \lceil 2 \sqrt{2N} \r \rceil  + 2} \nonumber \\
        &  =   &  \frac{2q^2+2}{\sqrt{2q}(q-1) + 2 \lceil 2\sqrt{2q^2+2} \rceil +4}. \nonumber
\end{eqnarray}

Next we consider the case $r=3$. Assume $q \ge 4$. Let $M$ be a subfield
of $L$ of degree $q/4$ over $K$.
Then this extension  is totally ramified at $\infty$ and
totally split at all other places of degree 1. Thus $M$ has $N=q^2/4+1$ rational 
places.  
Next we compute the genus $g(M)$ of $M$. From Theorem~2.1 of~\cite{gs4}
we have $g(M) = \sum_{i=1}^{t} E_i$ where $E_1$, $E_2$, ...,$E_t$ 
($t = q/4-1$) are the intermediate extensions $K \subseteq E_i \subseteq M$
with $[E_i:K]=2$. It follows from the proof of Proposition~1.2 of~\cite{gs4} 
that all the $E_i$ have genus  $g(E_i) = q_0$. Thus $g(M)=q_0(q/4 -1)$.

By~(\ref{brrr}) we have for  $q$ sufficiently large 
\begin{eqnarray}
B_3(M) 
&  \ge  &
\frac{q^3}{3}- \left (\frac{q}{q-1} + 2g(M) \frac{q^{1/2}}{q^{1/2}-1} \right ) 
\frac{q^{3/2}-1}{3} \nonumber \\
& \ge &
\frac{q^3}{3}  - \left ( 2 + 3g(M) \right ) \frac{ q^{3/2}-1}{3}  \nonumber \\
& = & 
\frac{q^3}{3}- \left ( 2 + 3q_0(q/4 -1) \right ) \frac{q^{3/2}-1}{3}\nonumber\\
& = & 
\frac{1}{3} \l (1-\frac{\sqrt{18}}{8} \r ) q^3 + \bigoh (q^{3/2}). \nonumber
\end{eqnarray}
Thus
$B_3 (M)>
\left\lfloor  3 + \lceil 2 \sqrt{2N} \rceil \right\rfloor
=\left\lfloor 3 + \lceil 2 \sqrt{q^2/2 + 2} \rceil  \right\rfloor$
 for all sufficiently large $q$.

As above, the ratio of class numbers $h(M \Fqr) /h(M)$ is odd.
Thus by Theorem~\ref{qeven}, for all $q$ sufficiently 
large, we get
\begin{eqnarray}
A(q^3)  & \ge  & \frac{N}{g(M) + \l \lceil 2 \sqrt{2N} \r \rceil  + 2}   \nonumber \\
        &  =   &  \frac{2q^2+8}{\sqrt{2q}(q-4) + 8 \lceil \sqrt{2q^2+8}\rceil +16},  \nonumber
\end{eqnarray}
as required.\endproof

\noindent {\em Remark}: 
The same ideas involved in the proof of the lower bound of $A(q^3)$ for 
$q$ even can be used to prove the following bounds  which improves the bounds of 
Corollary~\ref{cubes} and the bound~(\ref{cube5}) for 
characteristics 3, 5, and 7. 
\begin{thm}\label{357}
Let $q$ be a power of $p=$ $3$, $5$ or $7$. Then for all $q$ sufficiently large we have
\[    
A(q^3)   \ge  \frac{2(q^2+p^2)}{\sqrt{pq}(q-p^2) + 4p(p-1) \lceil \sqrt{q^2/p + p}\rceil +10p^2-12p}.
\]
\end{thm}
%
%


\subsection{Improvements of Serre's bound}\label{kemprf2}

\indent Throughout this section all logarithms will
be of base 2. First we assume that $q$ is odd. Let $r>0$ be an odd integer.
Put $k=\Fq (x)$, $k_r = \Fqr (x)$. Given $0 < \theta < 1/2$, let $n$ be the largest odd integer 
which does not exceed
$1+ \theta r \log q $. 
We choose $n$ to be odd merely for the sake of a neater proof.
Let $N_t = B_t(k)$ denote the number of monic irreducible polynomials 
of degree $t$ over $\Fq$. 
Let $m$ be the smallest integer such that 
$N_m \ge n$.

The lemma below shows that $m \le  \lceil 2\log r/ \log q \rceil +1$ 
for $q^r$ sufficiently large. 

\begin{lem}\label{pre1}
If $M= \lceil 2\log r/ \log q \rceil +1$, then $n \le N_M$ for 
all $q^r$ sufficiently large.
\end{lem}
\proof
Applying Proposition~\ref{brrrprop} to $F=k$,  we get
$N_M >(q^M - 2 q^{M/2})/M$.
Now
$q^M - 2 q^{M/2} \ge q^{2\log r/ \log q+1}-2q^{\log r/ \log q+1}
= qr(r -2)$. Since $qr(r -2)/M \ge 2 \frac{q(r-2)}{\log r+\log q}
\cdot r \log q$, the desired  result follows.\endproof

As $n \le N_m$, we may choose $n$ distinct monic irreducible  
polynomials
$P_1(x)$, $P_2(x)$, $\ldots$, $P_n(x) $ of degree $m$  over $\Fq$.
For $1\le i \le n$ define the extensions $k(y_i)/k$ with $y_i^2 = P_i(x)$. Let $H$ be the compositum of the fields $k(y_1), ..., k(y_n)$.
Further define the extension $k(y)$ by 
$y^2 = P_1(x)P_2(x) \cdots  P_n(x)$.
It is clear that the extension $H$ is 
an unramified abelian extension of $k(y)$ of exponent 2 and the Galois group has 2-rank equal to $n-1$.
Note that our choice of $n$ being odd ensures that  the place $\infty$  of $k$ 
does not ramify in the extension $H/k(y)$. 
By Proposition~\ref{brrrprop}, the number of degree $r$ places $B_r$ of $H$ 
satisfies
\begin{equation}\label{brrr2}
B_r > \frac{q^r}{r} - (2+7g(H)) \frac{q^{r/2}}{r}, 
\end{equation}
where $g(H)$ is the genus of $H$. Using the Hurwitz genus formula, we get $g(k(y))-1 = (mn+\epsilon-4)/2$ and $g(H)= 2^{n-2}(mn+\epsilon-4)+1$,
where $\epsilon$ is  1 if $m$ is odd and 0 otherwise.

Now let $P'$ be a place of degree $r$ in $H$ and let $P$ be the  place of
$k(y)$ which lies below $P'$. 
 Then $r=\deg P' = f(P'|P) \deg P$, where $f(P'|P)$ is the order of the 
decomposition group  $G(P'|P)$, which is cyclic of order at most 2.
Since $r$ is odd, we have $f(P'|P)=1$. Consequently the place $P$ splits
completely in the extension $H/k(y)$ and $\deg P = r$. 
Thus each degree $r$ place of $H$ divides a degree $r$ place of $k(y)$
which splits completely in $H/k(y)$.
Consequently the number of degree $r$
places of $k(y)$ which split completely in $H/k(y)$ is  $B_r/[H:k(y)] = B_r/2^{n-1}$.
From~(\ref{brrr2}) we have 
\[
B_r/2^{n-1} > \frac{q^r}{2^{n-1}r} - \l (\frac{9}{2^{n-1}}+\frac{7}{2}(mn-3) \r ) \frac{q^{r/2}}{r}. 
\]
As it is easily checked that
\[
\frac{q^r}{2^{n-1}r} - \l (\frac{9}{2^{n-1}}+\frac{7}{2}(mn-3) \r ) \frac{q^{r/2}}{r}
\ge \left(\frac{n-3}{2} \right )^2 -1 
\]
for all sufficiently large $q^r$, we can choose a set $S'$ of 
$ ( (n-3)/2 )^2  -1 $ places of degree $r$ of $k(y)$ which 
split completely in $H/k(y)$. Since $d_2 Cl_{S'} \ge n-1 =  2 + 2 \sqrt{|S'|+1}$,
we have by Proposition~\ref{golshaf} that $k(y)$ has an infinite $S'$-Hilbert class
field tower.
Passing to the constant field extension $k_r(y)$, each place of degree
$r$ in $k(y)$ splits into $r$ places of degree 1 in $k_r(y)$  and it is  easily seen that
$k_r(y)$ has an infinite $S$-Hilbert class field tower, where $S$ is the set 
of those places
of $k_r(y)$ which lie above those in $S'$. We thus get, again by 
Proposition~\ref{golshaf},
that
\begin{eqnarray}
A(q^r)& \ge & \frac{|S|}{g(k_r(y))-1} = \frac{|S'|r}{g(k(y))-1} =
2 (  (n-3)/2  )^2 -1 )r/(mn+\epsilon-4) \nonumber\\ 
& \ge &  ((n-3)^2 - 4)r/(2mn-6) \nonumber \\
& \ge &   
\frac{((\lfloor \theta r   \log q \rfloor-3)^2 - 4)r}
{2(\lceil 2\log r/ \log q \rceil +1)(\lfloor \theta r   \log q \rfloor +1)-6}\nonumber 
\end{eqnarray}
for all sufficiently large $q^r$.
We have proved
\begin{thm}\label{kemqodd}
Let $0 <  \theta < 1/2$.
Then for all sufficiently large odd $q^r$  we have 
\begin{equation}\label{obnd}
A(q^r) \ge
\frac{((\lfloor \theta r   \log q \rfloor-3)^2 - 4)r}
{2(\lceil 2\log r/ \log q \rceil +1)(\lfloor \theta r   \log q \rfloor +1)-6}.
\end{equation}
\end{thm}

For $q$ even the proof is essentially the same so we omit the details. 
The extensions $k(y_i)/k$ in this case are Artin-Schreier extensions 
defined by $y_i^2 +y_i = 1/P_i(x)$ and the extension
$k(y)/k$ is defined by $y^2 +y = \sum_i 1/P_i(x)$. Also in this case, $n$ 
need not be odd.
The lower bound we get in this case is approximately half that 
of the $q$ odd case. 
\begin{thm}\label{kemqeven}
Let $0 < \theta < 1/2$. Then 
for all sufficiently large  even $q^r$ we have 
\begin{equation}\label{ebnd}
A(q^r) \ge
\frac{\l ( \lfloor \theta r \log q \rfloor-2 \r )^2r}
{4(\lceil 2\log r/ \log q \rceil +1)(\lfloor \theta r \log q \rfloor +1)-8}.
\end{equation}
\end{thm}

As the right hand side of (\ref{obnd}) is at least
$\frac{r \log q}{2} \cdot 
\frac{\theta^2 r^2 (\log q)^2 - 8 \theta r \log q + 12 }
{\theta r \log q( \log r+ \log q) + \log r -2 \log q}$, which in turn is at least
$\frac{\theta }{4}r^2 (\log q)^2/(\log r + \log q)$ for all sufficiently 
large $q^r$,
we see that the bound~(\ref{obnd}) implies the bound~(\ref{kembd}).
The same is true for even $q$.

\section{ Lower bounds of $A(p)$ for small primes $p$} \label{smallp}

In view of the condition~(\ref{gscond}) in Proposition~\ref{golshaf}, it is 
important that we have 
good lower bounds for the $l$-rank of the $S$-divisor class group $Cl_S$.
Niederreiter and Xing~\cite{nx1} proved the following lower bound.
\begin{prop} \label{nxlb}
Let $F$ be a global function field and $K/F$  a finite abelian
extension. Let $T$ be a finite nonempty set of places of $F$ and $S$ the set of
places of $K$ lying over those in $T$. Then for any prime $l$ we have
\[
d_{l}Cl_{S}\ge\sum_{P}d_{l}G_{P}-d_{l}O_{T}^{*}-d_{l}G,
\]
 where  $G= \gal (K/F)$  and $G_P$  is the inertia subgroup at the
place $P$ of $F$. The sum is extended over all places $P$ of $F$.
\end{prop} 
Their proof uses Tate cohomology. Here we give another
 proof of this result assuming that at least one of 
the places in the set $T$ splits completely, which is the case in applications. 
The proof below, which uses
narrow ray class fields,  reveals that the lower bound of 
Proposition~\ref{nxlb}
is really a lower bound of the $l$-rank of the Galois group of the maximal subfield of the $S$-Hilbert class field of $K$ which is an abelian extension of $F$. If we remove the condition that a place of $T$ splits completely in $K/F$, then 
the proof below can be easily modified to obtain a lower bound which is
one less.
\begin{prop} \label{nxlb2}
Let $F/ \Fq$ be a global function field and $K/F$  a finite abelian
extension. Let $T$ be a finite nonempty set of places of $F$ and $S$ the set of
places of $K$ lying over those in $T$. 
If at least one place in 
$T$ splits completely in $K$, then  for any prime $l$ we have
\[
d_{l}Cl_{S}\ge 
       \sum_{P}d_{l}G_{P}-(|T|-1+d_{l}\mathbb{F}_{q}^*   )-d_{l}G,
\]
 where  $G= \gal (K/F)$, $G_P$  is the inertia subgroup at the
place $P$ of $F$. 
The sum is extended over all places $P$ of $F$.
\end{prop} 
{\em Remark}:  Observe that 
$d_{l}O_{T}^{*}= |T|-1+d_{l}\mathbb{F}_{q}^*$ so that  the bound coincides
with the one of Proposition~\ref{nxlb}. 

\noindent
{\em Proof of Proposition~\ref{nxlb2}}:
We continue with the notation introduced in section~\ref{nrcf}. 
Obviously, we may assume that the extension $K/F$ is ramified.
Denote by $\infty$ a place in $T$ which splits completely in $K$. 
Write $A$ for the ring of elements in $F$ regular outside $\infty$.
Let $M$ be the conductor of the extension $K/F$. Then 
$M$ is the smallest modulus for which $K$ is contained in the
narrow ray class field $F_M$.

Since all the field extensions involved are abelian, 
we may speak of the decomposition group or inertia group of 
places in the base field without specifying a corresponding place  
above. Let  $G_P''$ be the inertia group of a place $P$ 
of $F$ in the extension  $F_M/F$.
Now for any place $P'$ of $K$ which lies above a place $P$ in $F$, 
the inertia group of $P'$ in the extension $F_M/K$ is $G_P'' \cap \gal (F_M/K)$,
which is independent of the choice of $P'$. We denote the group 
$G_P'' \cap \gal (F_M/K)$ by $G_P'$. Observe that $G_\infty''$ is contained in $\gal(F_M/H_A)$ and $\gal(F_M/K)$. In particular, $G_{\infty}' = G_\infty''$.

If $J$ is the fixed field of $G_P''$, then
$G_P$ is isomorphic to 
$\gal (F_M / J \cap K)/ \gal (F_M/K)$$= [G_P'' \gal(F_M / K)] / \gal (F_M / K)$, which  is isomorphic to 
$G_P''/ G_P'' \cap \gal (F_M / K)$ 
$= G_P'' / G_P'$. In other words, $G_P \cong G_P'' / G_P'$.

Suppose that $M$ has prime decomposition $M=P_1^{e_1}P_2^{e_2} \cdots P_t^{e_t}$, where 
$P_1$, $P_2$, ..., $P_t$ are prime ideals of $A$ and 
$e_1, e_2, \, \ldots, e_t \ge 1$.
Let
$G'=G_{P_1}' \cdots G_{P_t}'$ and $G''=G_{P_1}'' \cdots G_{P_t}''$.
 Then $G'' = \gal (F_M/H_A)$ (cf. ~\cite{hay}, ~\cite{nx1}) so that $G_\infty'' \subseteq G''$. Let $L$ be the fixed
field of $G'G_\infty'$ in the extension $F_M/F$.
We have
\begin{eqnarray}\label{lbg}
d_{l} \gal (L/F) &=&  d_{l} \Pic_M(A)/ \gal (F_M/L) \nonumber \\
& \ge &  d_{l} G''/ G_\infty' G' \nonumber \\
 & = &  d_{l} \left ( G''/G' \right )/ \left ( G_\infty' G'/G' \right )
\nonumber \\
 &\ge & d_{l}  G''/G' - d_{l} G_\infty' G'/G' \nonumber \\
& = & \sum_{i=1}^t  d_{l} G_{P_i}-d_{l}  G_\infty' /  G_\infty' \cap G'
\nonumber \\
& =  & \sum_{P}  d_{l} G_{P}-d_{l}  G_\infty' /  G_\infty' \cap G'
,
\end{eqnarray}
where $P$ ranges over all places of $F$.

Since $\infty$ splits completely in $K$, the field $L$ contains $K$. Observe further that $L/K$ is an unramified abelian extension.
For each place $P \in T$ other than $\infty$, let $H_P'$ be the decomposition group of any place of $L$ dividing $P$, and let $H_P$ be the intersection of $H_P'$ with $\gal(L/K)$.
Denote by $K'$ the fixed field of all $H_P$, $P \in T$ and $P \ne \infty$. Then $K'$ is an unramified abelian extension of $K$ in which all places in $S$ split completely. In other words, $K'$ is a subfield of the $S$-Hilbert class field of $K$. Since each $H_P$ is cyclic, we have
\begin{equation}\label{ubg}
d_{l}\gal(L/K') \le |T| - 1.
\end{equation}

We now have 
\begin{eqnarray}
d_{l}Cl_{S} & \ge&  d_{l} \gal (K'/K)  =
             d_{l}  \gal (L/K)/ \gal (L/K')   \nonumber \\
           & \ge & d_{l}\gal (L/K) -  d_{l} \gal (L/K')  \nonumber \\
           & \ge  & d_{l}\gal (L/K)  - ( |T| - 1) \; \; (\hbox{by } ~(\ref{ubg}))\nonumber \\
           & \ge  & d_{l}\gal(L/F) - d_{l}\gal (K/F)  - ( |T|-1) \nonumber \\
           & \ge  & 
\sum_{P}  d_{l} G_{P}-( |T|-1 +d_{l}  G_\infty' /  G_\infty' \cap G') - d_{l}G 
 \; \; (\hbox{by } ~(\ref{lbg})) \nonumber \\
 & \ge &  \sum_{P}  d_{l} G_{P}-(|T|-1 +d_{l}  G_\infty' ) - d_{l}G.
\nonumber
\end{eqnarray}
Since $ G_\infty' \cong \mathbb{F}_{q}^*$, we are done.\endproof

\noindent
Next we present lower bounds for $A(p)$ where $p=7,11,13, 17$.

\begin{thm}
We have 
\[ A(7) \ge 9/10 \]
\end{thm}
\proof
Let $k$ be the rational function field $\Fqseven (x)$. Let $F=k(y)$ be the 
function field defined by
\[
y^2=Q(x):=x^6+2x^5+3x^4+3x^3+x^2+1.
\]
Then $F/k$ is a Kummer extension in which all the rational  places of 
$k$  split completely.
The only place ramifying
in $F/k$ is  $Q(x)$ and by the Hurwitz genus formula 
$g(F)=2$.
  
Let $K=k(z)$ be the function field defined by $z^2= P(x)$ where
\begin{eqnarray} \nonumber
P(x) & = &  x(x+1)(x+2)(x^2+4x+6)(x^2+3x+6)(x^2+3x+1) \\ \nonumber
     &   & (x^2+6x+4)(x^2+6x+3)(x^2+2x+2)(x^2+4). \nonumber
\end{eqnarray}
Then $K/k$ is a Kummer extension in which the places $x+3, x+4, x+5, x+6$ 
split completely and the  ramified places are those in 
the set $R = \{ 
x,x+1 , x+2 , x^2+4x+6 , x^2+3x+6 , x^2+3x+1 , x^2+6x+4 , x^2+6x+3 , x^2+2x+2  , x^2+4 , \infty \}$.
Now, from the  relations 
\begin{eqnarray}
Q(x) & \equiv   & (3+2x)^2    \hbox{ mod }  x^2+4x+6  \nonumber \\
Q(x) & \equiv   & (3+x)^2   \hbox{ mod }  x^2+3x+6  \nonumber \\
Q(x) & \equiv   & (2+5x)^2    \hbox{ mod }  x^2+3x+1  \nonumber \\
Q(x) & \equiv   & 1^2    \hbox{ mod }  x^2+6x+4  \nonumber \\
Q(x) & \equiv   & (2+2x)^2    \hbox{ mod }  x^2+6x+3  \nonumber \\
Q(x) & \equiv   & (3+x)^2    \hbox{ mod }  x^2+2x+2  \nonumber \\
Q(x) & \equiv   & (2+5x)^2    \hbox{ mod }  x^2+4  \nonumber 
\end{eqnarray} 
it follows that all the places in the set $R$ split completely
in the extension $F/k$.
Let $T$ be the set of places of $F$ lying over
$x+2, x+3, x+4, x+5, x+6$ and $S$ the set of places of $FK$ lying over 
those in $T$. Then $|T|=2\cdot5= 10$ and
$|S|=2 + 2 \cdot 8= 18$. By Proposition~\ref{nxlb},
$$d_{2}Cl_{S}\ge\sum_{P}d_{2}G_{P}-|T|-d_{2}G=22-10-1=11,$$
where $G=$ Gal$(FK/F)\simeq{\bf Z}/2{\bf Z}$ and the sum runs over
all places $P$ of $F$. Since $11 \ge 2+2\sqrt{|S| + 1}$, 
the condition~(\ref{gscond}) in Proposition~\ref{golshaf} is satisfied.
By the Hurwitz genus formula we have
$2g(FK)-2=2(2g(F)-2)+ 2(4\cdot 1 + 7 \cdot 2) =40$,
and so
$$A(7)\ge\frac{|S|}{g(FK)-1}=\frac{9}{10}.\Box$$\medskip

\begin{thm}
We have 
\[ A(11) \ge 12/11 = 1.0909 \ldots \]
\end{thm}
\proof
Put $k=\FqXI (x)$. 
Let \newline
$P(x) = 
(x^2+4x+2)(x^2+5x+7)(x^2+8x+9)(x^2+6x+7)(x^2+1)(x^2+3)(x^2+4)(x^2+5)
(x^2+9)(x^2+10x+6)(x^2+6x+3)(x^2+x+1)(x^2+6x+2)(x^2+9x+5)(x^2+6x+10)
(x^2+x+4)(x^2+x+6)(x^2+x+7)(x^2+x+8)(x^2+10x+4)(x^2+9x+4)(x^2+9x+10)
(x^2+6x+1)(x^2+7x+9)$,
which is a product of 24  irreducible polynomials of degree 2  over $\FqXI$,
call them $P_1 (x)$,..., $P_{24}(x)$.

Consider the extension $k(y)$ defined by $y^2=P(x)$. 
Now $k(y)$ is contained in the function field $F=k(y_1 , \ldots ,y_{24})$,
where $y_i^2 = P_i(x)$ for $ 1 \le i \le 24$.
Moreover the extension $F/k(y)$ is unramified, $\gal (F/k(y)) \cong 
(\ZZ / 2\ZZ )^{23}$ and the place $\infty$ splits completely in 
$F/k$.

Now the the places
in the set $T= \{ x+ \alpha | \alpha \in \FqXI \} \cup \{ \infty \}$
split completely in $k(y)/k$. Thus $k(y)$ is contained in the
 decomposition fields of the places in $T$ .
For each place $x+ \alpha$ in $T$ let 
$G_\alpha$ be the decomposition group of $x+ \alpha$ in $\gal (K/k)$.
Then  $G_\alpha$ is a cyclic subgroup of $\gal (F/k(y))$.

Let  $H$ be the the subgroup of $\gal (F/k(y))$
generated by the groups $G_\alpha$.  
Since each group $G_\alpha$  is  cyclic of order at most  2, it follows that
$d_{2} H \le 11$.
Let $K'$ be the fixed field of $H$ in $F$ and let $S$ be the set of 
places in $k(y)$ which lie above those in $T$. 
Then $K'/k(y)$ is an unramified abelian extension in which each place in 
$S$ splits completely.
We now have $d_{2}Cl_{S} \ge d_{2}\gal (K'/k(y)) \ge d_{2}\gal (F/k(y))-
d_{2}\gal (F/K') \ge 12 = 2 +2\sqrt{|S|+1} $.

By the Hurwitz genus formula $g(k(y))-1 = \frac{1}{2}(-4 + 24 \cdot2 )   =22$.
Thus  by Proposition~\ref{golshaf} we have $A(11) \ge  24/22 = 1.0909 \ldots $.\endproof
\begin{thm}
We have
\[
A(13)  \ge 4/3 = 1.333 \ldots
\]
\end{thm}
\proof
Put $k = \FqXIII (x)$ and let 
 $ P(x) = x(x-1)(x-2)(x-3)(x-4)(x-5)(x-6)(x-7)(x-9) $.
Define the extension $k(y)/k$ by $y^2=P(x)$. Then $g(k(y)) = 3$ and the only rational places that
split completely in $k(y)/k$ are those in the set $T= \{x+2,  x+3\}$. If $S$ is the set of 
4 places in $k(y)$ which lie above those in $T$, then by Proposition~\ref{nxlb}, we have $d_2 Cl_S \ge 10 - 2- 1 =7$. Since $7 > 2 + 2\sqrt{|S|+1}$, we
have from Proposition~\ref{golshaf} that 
$A(13)  \ge |S|/(g(k(y))-1) = 4/3$ as required.\endproof

Likewise, by using the polynomial 
$P(x) = x(x-1)(x-2)(x-3)(x-4)(x-5)(x-6)(x-7)(x-8)(x-9)(x-11)(x-12)(x-15)$ one can show
that $A(17) \ge 8/5$.

\section{Acknowledgements}

We  thank H. Niederreiter and A. Temkine for providing us  with preprints of their papers related to 
the topic of this paper.


\end{document}